\documentclass{amsart}
\usepackage{setspace}
\usepackage{a4}
\usepackage{amssymb,amsmath,amsthm,latexsym}
\usepackage{amsfonts}
\usepackage{amsfonts}
\usepackage{graphicx}
\usepackage{textcomp}
\usepackage{cite}
\usepackage{enumerate}
\usepackage[mathscr]{euscript}
\usepackage{mathtools}
\newtheorem{theorem}{Theorem}[section]

\newtheorem{conjecture}[theorem]{Conjecture}

\newtheorem{problem}[theorem]{Problem}

\newtheorem{remark}[theorem]{Remark}

\title{This is the title}
\usepackage{hyperref}
\usepackage{enumerate}
\usepackage{mathtools}
\usepackage{amsmath}
\usepackage{tikz}
\usepackage{hyperref}
\usepackage{graphicx}
\usepackage{textcomp}
\usetikzlibrary{cd}
\usepackage{fancyhdr}

\begin{document}
\begin{center}
{\bf{MODULAR BOURGAIN-TZAFRIRI RESTRICTED INVERTIBILITY CONJECTURES AND JOHNSON-LINDENSTRAUSS FLATTENING   CONJECTURE}}\\
\vspace{0.2cm}
\hrule 
\vspace{0.2cm}
\textbf{K. MAHESH KRISHNA} \\
Post Doctoral Fellow \\
Statistics and Mathematics Unit\\
Indian Statistical Institute, Bangalore Centre\\
Karnataka 560 059 India\\
Email: kmaheshak@gmail.com \\
Date: \today
\end{center}

\hrule
\vspace{0.5cm}
%--------------------------------------
\textbf{Abstract}:  We recently formulated  important Modular Bourgain-Tzafriri Restricted Invertibility Conjectures and Modular Johnson-Lindenstrauss Flattening  Conjecture in the Appendix of  \textit{[arXiv: 2207.12799.v1]}. For the sake of wide accessibility we give a self-contained treatment of them.


\hrule
\tableofcontents
\hrule

\section{Modular Bourgain-Tzafriri Restricted Invertibility Conjectures}
 Let $d\in \mathbb{N}$, $ \mathbb{K}=\mathbb{C}$ or $\mathbb{R}$ and $\mathbb{K}^d$ be the standard $d$-dimensional Hilbert space with canonical orthonormal basis $\{e_j\}_{j=1}^d$.  Following result obtained by  Prof. Bourgain and Prof. Tzafriri  in 1987 is one of the highly useful  result of \textbf{20 century Mathematics} which combines algebraic property of a matrix  (invertibility) with the analytic property (norm) \cite{BOURGAINTZAFRIRI}.
\begin{theorem}\cite{BOURGAINTZAFRIRI, BOURGAINTZAFRIRI2}\label{BOURGAINTZARIRITHEOREM}  \textbf{(Bourgain-Tzafriri Restricted Invertibility Theorem)
		There are universal constants $A>0$, $c>0$ satisfying the following property. If  $d \in \mathbb{N}$, and  $T:\mathbb{K}^d \to \mathbb{K}^d$ is a linear operator   with  $\|Te_j\|=1$, $\forall 1\leq j\leq d$, then there exists a subset $\sigma \subseteq \{1, \dots, d\}$ of cardinality 
		\begin{align*}
			\operatorname{Card}(\sigma)\geq \frac{cd}{\|T\|^2}
		\end{align*}
		such that 
		\begin{align*}
			\left\|\sum_{j \in \sigma}a_jTe_j\right\|^2\geq A \sum_{j \in \sigma}|a_j|^2, \quad \forall a_j \in  \mathbb{K}, \forall j \in \sigma.
	\end{align*}}
\end{theorem}  
In 2009 Tropp gave a polynomial time algorithm for the proof Theorem \ref{BOURGAINTZARIRITHEOREM} \cite{TROPP}. It came as a surprise in 2009 (arXiv version) when  Spielman and Srivastava gave a simple proof of Theorem \ref{BOURGAINTZARIRITHEOREM} by proving a generalization of  Theorem \ref{BOURGAINTZARIRITHEOREM}  \cite{SPIELMANSRIVASTAVA} due to Vershynin \cite{VERSHYNIN}.  In 2012 Casazza and Pfander  proved infinite dimensional version of Theorem \ref{BOURGAINTZARIRITHEOREM} \cite{CASAZZAPFANDER}. For  beautiful descriptions of Bourgain-Tzafriri restricted invertibility theorem and generalizations we refer \cite{NAORYOUSSEF, CASAZZATREMAIN, NAOR, SRIVASTAVA, CASAZZA, YOUSSEF, RAVICHANDRAN, TROPP2, CHRETIENDARSES, YOUSSEF2, VERSHYNIN, MARCUSSPIELMANSRIVASTAVA, GLUSKINOLEVSKII}.

We formulate Conjectures \ref{RC} and \ref{NRC} which are  based on Theorem \ref{BOURGAINTZARIRITHEOREM}. Let $d\in \mathbb{N}$, $\mathcal{A}$ be a  unital C*-algebra and $\mathcal{A}^d$ be the left module over $\mathcal{A}$ w.r.t. natural operations. Modular $\mathcal{A}$-inner product on $\mathcal{A}^d$ is defined as 
  \begin{align*}
  	\langle (x_j)_{j=1}^d, (y_j)_{j=1}^d\rangle \coloneqq \sum_{j=1}^{d}x_jy_j^*,\quad \forall (x_j)_{j=1}^d, (y_j)_{j=1}^d \in \mathcal{A}^d.
  \end{align*}
  Hence the norm on $\mathcal{A}^d$ becomes 
  \begin{align*}
  	\|(x_j)_{j=1}^d\|\coloneqq \left\|\sum_{j=1}^{d}x_jx_j^*\right\|^\frac{1}{2}, \quad \forall (x_j)_{j=1}^d \in \mathcal{A}^d.
  \end{align*}
In this way $\mathcal{A}^d$ becomes \textbf{standard Hilbert C*-module} \cite{KAPLANSKY, PASCHKE, RIEFFEL}.  
  \begin{conjecture}\label{RC} 
  \textbf{[(Commutative) Modular Bourgain-Tzafriri Restricted Invertibility Conjecture]
  Let $\mathcal{A}$ be a  unital commutative C*-algebra and $\mathcal{I}(\mathcal{A})$ be the set of all invertible elements of $\mathcal{A}$.   For $d\in \mathbb{N}$, let $ \mathbb{M}_{d\times d}(\mathcal{A})$ be the set of all $d$ by $d$ matrices over $\mathcal{A}$. For  $ M\in \mathbb{M}_{d\times d}(\mathcal{A})$, let $\det (M) $ be the determinant of $M$. Let $\mathcal{A}^d$ be the standard (left) Hilbert C*-module over 	$\mathcal{A}$ and    $\{e_j\}_{j=1}^d$ be the  canonical orthonormal basis for $\mathcal{A}^d$. There are universal real  constants $A>0$, $c>0$ ($A$ and $c$ may  depend upon C*-algebra $\mathcal{A}$) satisfying the following property. If  $d \in \mathbb{N}$ and  $ M \in \mathbb{M}_{d\times d}(\mathcal{A})$    with  $\langle Me_j, Me_j \rangle =1$, $\forall 1\leq j\leq d$ and $\det (M) \in \mathcal{I}(\mathcal{A})\cup\{0\}$, then there exists a subset $\sigma \subseteq \{1, \dots, d\}$ of cardinality 
  \begin{align*}
  	\operatorname{Card}(\sigma)\geq \frac{cd}{\|M\|^2}
  \end{align*}
  such that 
  \begin{align*}
  	\sum_{j \in \sigma}\sum_{k \in \sigma}a_j \langle Me_j, Me_k\rangle a_k^*=\left\langle \sum_{j \in \sigma}a_jMe_j, \sum_{k \in \sigma}a_kMe_k\right\rangle \geq A \sum_{j \in \sigma}a_ja_j^*, \quad \forall a_j \in  \mathcal {A}, \forall j \in \sigma, 
  \end{align*}	
  where $\|M\|$ is the norm of the Hilbert C*-module homomorphism defined by $M$ as $M:\mathcal{A}^d \ni x \mapsto Mx \in \mathcal{A}^d$.}
   \end{conjecture}
  We next formulate Conjecture \ref{RC} for unital C*-algebras which need not be commutative. In the statement we use the notion of Manin matrices. We refer \cite{CHERVOVFALQUIRUBTSOV, CHERVOVFALQUI} for the basics of Manin matrices.
  \begin{conjecture}\label{NRC}
   \textbf{[(Noncommutative) Modular Bourgain-Tzafriri Restricted Invertibility Conjecture ]
  	Let $\mathcal{A}$ be a  unital  C*-algebra and $\mathcal{I}(\mathcal{A})$ be the set of all invertible elements of $\mathcal{A}$.   For $d\in \mathbb{N}$, let $ \mathbb{MM}_{d\times d}(\mathcal{A})$ be the set of all $d$ by $d$ Manin matrices over $\mathcal{A}$. For  $ M\in \mathbb{MM}_{d\times d}(\mathcal{A})$, let $\det^\text{column} (M) $ be the Manin determinant of $M$ by column expansion. Let $\mathcal{A}^d$ be the standard (left) Hilbert C*-module over 	$\mathcal{A}$ and    $\{e_j\}_{j=1}^d$ be the  canonical orthonormal basis for $\mathcal{A}^d$. There are universal real  constants $A>0$, $c>0$ ($A$ and $c$ may  depend upon C*-algebra $\mathcal{A}$) satisfying the following property. If  $d \in \mathbb{N}$ and  $ M \in \mathbb{MM}_{d\times d}(\mathcal{A})$    with  $\langle Me_j, Me_j \rangle =1$, $\forall 1\leq j\leq d$ and $\det^\text{column} (M) \in \mathcal{I}(\mathcal{A})\cup\{0\}$, then there exists a subset $\sigma \subseteq \{1, \dots, d\}$ of cardinality 
  	\begin{align*}
  		\operatorname{Card}(\sigma)\geq \frac{cd}{\|M\|^2}
  	\end{align*}
  	such that 
  	\begin{align*}
  		\sum_{j \in \sigma}\sum_{k \in \sigma}a_j \langle Me_j, Me_k\rangle a_k^*=\left\langle \sum_{j \in \sigma}a_jMe_j, \sum_{k \in \sigma}a_kMe_k\right\rangle \geq A \sum_{j \in \sigma}a_ja_j^*, \quad \forall a_j \in  \mathcal {A}, \forall j \in \sigma, 
  	\end{align*}	
  	where $\|M\|$ is the norm of the Hilbert C*-module homomorphism defined by $M$ as $M:\mathcal{A}^d \ni x \mapsto Mx \in \mathcal{A}^d$.} 
  \end{conjecture}
 \begin{remark}
 	\begin{enumerate}[\upshape (i)]
 		\item \textbf{We can surely formulate Conjecture \ref{RC} by removing the condition $\det (M) \in \mathcal{I}(\mathcal{A})\cup\{0\}$ and Conjecture \ref{NRC} by removing the condition Manin matrices and $\det^\text{column} (M) \in \mathcal{I}(\mathcal{A})\cup\{0\}$. But we strongly  believe that Conjectures  \ref{RC} and \ref{NRC}  fail with this much of generality.}
 		\item \textbf{If Conjecture \ref{RC} holds but Conjecture \ref{NRC} fails,  then we  can try Conjecture \ref{NRC} for W*-algebras or C*-algebras with invariant basis number (IBN) property.} We refer \cite{GIPSON} for the IBN properties of C*-algebras.
 	\end{enumerate}
 \end{remark}

  \section{Modular Johnson-Lindenstrauss Flattening Conjecture}
   Everything starts from the following surprising result of Prof. Johnson and Prof. Lindenstrauss from 1984 \cite{JOHNSONLINDENSTRAUSS}.
  \begin{theorem}\cite{MATOUSEKBOOK, JOHNSONLINDENSTRAUSS}\label{JL}
  	\textbf{(Johnson-Lindenstrauss Flattening Lemma - original form) Let $M, N\in \mathbb{N}$ and  $\mathbf{x}_1, \mathbf{x}_2, \dots, \mathbf{x}_M\in \mathbb{R}^N$. For each $0<\varepsilon<1$, there exists a Lipschitz map $f: \mathbb{R}^N\to \mathbb{R}^m$  and a real $r>0$ such that 
  		\begin{align*}
  			r(1-\varepsilon)\|\mathbf{x}_j-\mathbf{x}_k\|\leq \|f(\mathbf{x}_j)-f(\mathbf{x}_k)\|\leq r(1+\varepsilon)\|\mathbf{x}_j-\mathbf{x}_k\|, \quad \forall 1 \leq j, k\leq M, 
  		\end{align*}
  		where 
  		\begin{align*}
  			m=O\left(\frac{\log M}{\varepsilon^2}\right).
  	\end{align*}}
  \end{theorem}
Over the time, several improvements of Theorem \ref{JL} were obtained. We recall them. 
  \begin{theorem}\cite{FRANKLMAEHARA, JOHNSONLINDENSTRAUSS} \textbf{(Johnson-Lindenstrauss  Flattening  Lemma - Frankl-Maehara form)
  		Let $0<\varepsilon<\frac{1}{2}$   and $ M\in \mathbb{N}$. Define 
  		\begin{align*}
  			m (\varepsilon, M)\coloneqq \left\lceil 9\frac{1}{\varepsilon^2-\frac{2\varepsilon^3}{3}}\log M \right \rceil +1.
  		\end{align*} 
  		If  $M> m (\varepsilon, M)$, then for any  $\mathbf{x}_1, \mathbf{x}_2, \dots, \mathbf{x}_M\in \mathbb{R}^M$,  there exists 	a map $f: \{\mathbf{x}_j\}_{j=1}^M\to \mathbb{R}^m$   such that 
  		\begin{align*}
  			(1-\varepsilon)\|\mathbf{x}_j-\mathbf{x}_k\|^2\leq \|f(\mathbf{x}_j)-f(\mathbf{x}_k)\|^2\leq (1+\varepsilon)\|\mathbf{x}_j-\mathbf{x}_k\|^2, \quad \forall 1 \leq j, k\leq M.
  	\end{align*}}
  \end{theorem}
  \begin{theorem}\cite{DASGUPTAGUPTA, JOHNSONLINDENSTRAUSS}
  	\textbf{(Johnson-Lindenstrauss Flattening Lemma - Dasgupta-Gupta form)	
  		Let $M, N\in \mathbb{N}$ and  $\mathbf{x}_1, \mathbf{x}_2, \dots, \mathbf{x}_M\in \mathbb{R}^N$. Let $0<\varepsilon<1$. Choose any natural number $m$ such that 
  		\begin{align*}
  			m>4\frac{1}{\frac{\varepsilon^2}{2}-\frac{\varepsilon^3}{3}}\log M. 
  		\end{align*}
  		Then there exists a  map $f: \mathbb{R}^N\to \mathbb{R}^m$   such that 
  		\begin{align*}
  			(1-\varepsilon)\|\mathbf{x}_j-\mathbf{x}_k\|^2\leq \|f(\mathbf{x}_j)-f(\mathbf{x}_k)\|^2\leq (1+\varepsilon)\|\mathbf{x}_j-\mathbf{x}_k\|^2, \quad \forall 1 \leq j, k\leq M.
  		\end{align*}
  		The map $f$ can be found in randomized polynomial time.}
  \end{theorem}
  \begin{theorem} \cite{JOHNSONLINDENSTRAUSS, FOUCARTRAUHUT} \textbf{(Johnson-Lindenstrauss Flattening Lemma - matrix form)
  		There is a universal constant $C>0$ satisfying the following. Let $0<\varepsilon<1$, $M,N \in \mathbb{N}$ and $\mathbf{x}_1, \mathbf{x}_2, \dots, \mathbf{x}_M\in \mathbb{R}^N$. For each natural number 
  		\begin{align*}
  			m>\frac{C}{\varepsilon^2}\log M,
  		\end{align*}
  		there exists a matrix $M\in \mathbb{M}_{m\times N}(\mathbb{R})$ such that 
  		\begin{align*}
  			(1-\varepsilon)\|\mathbf{x}_j-\mathbf{x}_k\|^2\leq \|M(\mathbf{x}_j-\mathbf{x}_k)\|^2\leq (1+\varepsilon)\|\mathbf{x}_j-\mathbf{x}_k\|^2, \quad \forall 1 \leq j, k\leq M.
  	\end{align*}}
  \end{theorem}
  In 2016 Larsen and Nelson\cite{LARSENNELSON, LARSENNELSON2017} showed that bound in Johnson-Lindenstrauss flattening lemma is optimal. For a comprehensive look on Johnson-Lindenstrauss flattening lemma from the historical point of view and its various applications we refer important papers and books  \cite{ACHLIOPTAS, INDYKMOTWANI, KRAHMERWARD, INDYKNAOR, ALON, COHENELDERMUSCOMUSCOPERSU, BECCHETTIBURYCOHENADDADGRANDONISCHWIEGELSHOHN, JAYRAMWOODRUFF, OSTROVSKII, LINIALLONDONRABINOVICH, DEDEZALAURENT, BOURGAIN, JACQUES, INDYK2000, AILONCHAZELLE, AILONCHAZELLE2009, KANENELSON, DASGUPTAKUMARSARLOS, BARTALRECHTSCHULMAN, KANEMEKANELSON, KLARTAGMENDELSON, ARRIAGAVEMPALA, AILONLIBERTY, AILONLIBERTY2, ALONKLARTAG, FREKSEN, VENKATASUBRAMANIANWANG, LIBERTYAILONSINGER}. \\
  We now  formulate the following interesting problem. 
  \begin{problem}\label{MODULARPROBLEM}
  		\textbf{Let $\mathscr{A}$ be the set of all unital C*-algebras. What is the best function $\phi: \mathscr{A}\times (0,1)\times \mathbb{N} \to (0,\infty)$ satisfying the following. Let $\mathcal{A}$ be a unital C*-algebra. 	There is a universal constant $C>0$ (which may depend upon $\mathcal{A}$) satisfying the following. Let $0<\varepsilon<1$, $M,N \in \mathbb{N}$ and $\mathbf{x}_1, \mathbf{x}_2, \dots, \mathbf{x}_M\in \mathcal{A}^N$. For each natural number 
  	\begin{align*}
  		m>C\phi (\mathscr{A}, \varepsilon, M)
  	\end{align*}
  	there exists a matrix $M\in \mathbb{M}_{m\times N}(\mathcal{A})$ such that 
  	\begin{align*}
  		(1-\varepsilon)\langle \mathbf{x}_j-\mathbf{x}_k, \mathbf{x}_j-\mathbf{x}_k\rangle \leq \langle M(\mathbf{x}_j-\mathbf{x}_k), M(\mathbf{x}_j-\mathbf{x}_k)\rangle \leq (1-\varepsilon)\langle \mathbf{x}_j-\mathbf{x}_k, \mathbf{x}_j-\mathbf{x}_k\rangle, \quad \forall 1 \leq j, k\leq M.
  \end{align*}}
  \end{problem}
   A particular case of Problem \ref{MODULARPROBLEM} is the following  conjecture.
 \begin{conjecture}\label{MJLC}
	\textbf{(Modular Johnson-Lindenstrauss Flattening Conjecture)
		Let $\mathcal{A}$ be a unital C*-algebra. 	There is a universal constant $C>0$ (which may depend upon $\mathcal{A}$) satisfying the following. Let $0<\varepsilon<1$, $M,N \in \mathbb{N}$ and $\mathbf{x}_1, \mathbf{x}_2, \dots, \mathbf{x}_M\in \mathcal{A}^N$. For each natural number 
		\begin{align*}
			m>\frac{C}{\varepsilon^2}\log M,
		\end{align*}
		there exists a matrix $M\in \mathbb{M}_{m\times N}(\mathcal{A})$ such that 
		\begin{align*}
			(1-\varepsilon)\langle \mathbf{x}_j-\mathbf{x}_k, \mathbf{x}_j-\mathbf{x}_k\rangle \leq \langle M(\mathbf{x}_j-\mathbf{x}_k), M(\mathbf{x}_j-\mathbf{x}_k)\rangle \leq (1-\varepsilon)\langle \mathbf{x}_j-\mathbf{x}_k, \mathbf{x}_j-\mathbf{x}_k\rangle, \quad \forall 1 \leq j, k\leq M.
	\end{align*}}	
\end{conjecture}  
\begin{remark}
	\begin{enumerate}[\upshape(i)]
		\item We believe that Conjecture \ref{MJLC}  holds  at least for W*-algebras (von Neumann algebras) or C*-algebras with IBN property. 	
		\item 	\textbf{Modular Welch bounds} are derived in \cite{MAHESHKRISHNA} and various problems including modular Zauner's conjecture are stated there. 
	\end{enumerate}
\end{remark}

 \bibliographystyle{plain}
 \bibliography{reference.bib}

\end{document}